\documentclass[12pt]{amsart}

\usepackage{amsmath,amssymb,amsthm,amscd,a4wide,upref}
\usepackage{mathrsfs}
\usepackage[enableskew]{youngtab}

\usepackage{hyperref}
\hypersetup{colorlinks,citecolor=blue,filecolor=black,linkcolor=blue,urlcolor=blue}
\usepackage{enumerate}
\usepackage{mathtools}
\usepackage{enumitem}
\usepackage[dvipsnames]{xcolor}

\usepackage{lmodern}
\usepackage[T1]{fontenc}

\begin{document}
	
\newcommand{\End}{{\rm{End}\ts}}
\newcommand{\Hom}{{\rm{Hom}}}
\newcommand{\Mat}{{\rm{Mat}}}
\newcommand{\ch}{{\rm{ch}\ts}}
\newcommand{\sh}{{\rm{sh}}}
\newcommand{\chara}{{\rm{char}\ts}}
\newcommand{\diag}{ {\rm diag}}
\newcommand{\non}{\nonumber}
\newcommand{\wt}{\widetilde}
\newcommand{\wh}{\widehat}
\newcommand{\ot}{\otimes}
\newcommand{\la}{\lambda}
\newcommand{\La}{\Lambda}
\newcommand{\De}{\Delta}
\newcommand{\al}{\alpha}
\newcommand{\be}{\beta}
\newcommand{\ga}{\gamma}
\newcommand{\Ga}{\Gamma}
\newcommand{\ep}{\epsilon}
\newcommand{\ka}{\kappa}
\newcommand{\vk}{\varkappa}
\newcommand{\si}{\sigma}
\newcommand{\vp}{\varphi}
\newcommand{\de}{\delta}
\newcommand{\ze}{\zeta}
\newcommand{\om}{\omega}
\newcommand{\ee}{\epsilon^{}}
\newcommand{\su}{s^{}}
\newcommand{\hra}{\hookrightarrow}
\newcommand{\ve}{\varepsilon}
\newcommand{\ts}{\,}
\newcommand{\tss}{\hspace{1pt}}
\newcommand{\vac}{\mathbf{1}}
\newcommand{\vacr}{|\tss 0\rangle}
\newcommand{\vacl}{\langle 0\tss |}
\newcommand{\di}{\partial}
\newcommand{\qin}{q^{-1}}
\newcommand{\Sr}{ {\rm S}}
\newcommand{\U}{ {\rm U}}
\newcommand{\BL}{ {\overline L}}
\newcommand{\BE}{ {\overline E}}
\newcommand{\BP}{ {\overline P}}
\newcommand{\AAb}{\mathbb{A}\tss}
\newcommand{\CC}{\mathbb{C}\tss}
\newcommand{\KK}{\mathbb{K}\tss}
\newcommand{\QQ}{\mathbb{Q}\tss}
\newcommand{\SSb}{\mathbb{S}\tss}
\newcommand{\ZZ}{\mathbb{Z}\tss}
\newcommand{\X}{ {\rm X}}
\newcommand{\Y}{ {\rm Y}}
\newcommand{\Z}{{\rm Z}}
\newcommand{\Ac}{\mathcal{A}}
\newcommand{\achi}{\Ac_{\chi}}
\newcommand{\bachi}{\overline\Ac_{\chi}}
\newcommand{\Lc}{\mathcal{L}}
\newcommand{\ol}{\overline}
\newcommand{\Pc}{\mathcal{P}}
\newcommand{\Qc}{\mathcal{Q}}
\newcommand{\Tc}{\mathcal{T}}
\newcommand{\Sc}{\mathcal{S}}
\newcommand{\Bc}{\mathcal{B}}
\newcommand{\Dc}{\mathcal{D}}
\newcommand{\Ec}{\mathcal{E}}
\newcommand{\Fc}{\mathcal{F}}
\newcommand{\Hc}{\mathcal{H}}
\newcommand{\Uc}{\mathcal{U}}
\newcommand{\Vc}{\mathcal{V}}
\newcommand{\Wc}{\mathcal{W}}
\newcommand{\Yc}{\mathcal{Y}}
\newcommand{\Ar}{{\rm A}}
\newcommand{\Br}{{\rm B}}
\newcommand{\Ir}{{\rm I}}
\newcommand{\Fr}{{\rm F}}
\newcommand{\Jr}{{\rm J}}
\newcommand{\Or}{{\rm O}}
\newcommand{\GL}{{\rm GL}}
\newcommand{\Spr}{{\rm Sp}}
\newcommand{\Rr}{{\rm R}}
\newcommand{\Zr}{{\rm Z}}
\newcommand{\g}{\mathfrak{g}}
\newcommand{\gl}{\mathfrak{gl}}
\newcommand{\middd}{{\rm mid}}
\newcommand{\ev}{{\rm ev}}
\newcommand{\Pf}{{\rm Pf}}
\newcommand{\Norm}{{\rm Norm\tss}}
\newcommand{\oa}{\mathfrak{o}}
\newcommand{\spa}{\mathfrak{sp}}
\newcommand{\osp}{\mathfrak{osp}}
\newcommand{\h}{\mathfrak h}
\newcommand{\n}{\mathfrak n}
\newcommand{\z}{\mathfrak{z}}
\newcommand{\Zgot}{\mathfrak{Z}}
\newcommand{\p}{\mathfrak{p}}
\newcommand{\sll}{\mathfrak{sl}}
\newcommand{\agot}{\mathfrak{a}}
\newcommand{\qdet}{ {\rm qdet}\ts}
\newcommand{\Ber}{ {\rm Ber}\ts}
\newcommand{\HC}{ {\mathcal HC}}
\newcommand{\cdet}{ {\rm cdet}}
\newcommand{\tr}{ {\rm tr}}
\newcommand{\gr}{ {\rm gr}}
\newcommand{\str}{ {\rm str}}
\newcommand{\st}{ {\rm st}}
\newcommand{\loc}{{\rm loc}}
\newcommand{\Gr}{{\rm G}}
\newcommand{\sgn}{ {\rm sgn}\ts}
\newcommand{\ba}{\bar{a}}
\newcommand{\bb}{\bar{b}}
\newcommand{\bi}{\bar{\imath}}
\newcommand{\bj}{\bar{\jmath}}
\newcommand{\bk}{\bar{k}}
\newcommand{\bl}{\bar{l}}
\newcommand{\hb}{\mathbf{h}}
\newcommand{\Sym}{\mathfrak S}
\newcommand{\fand}{\quad\text{and}\quad}
\newcommand{\Fand}{\qquad\text{and}\qquad}
\newcommand{\For}{\qquad\text{or}\qquad}
\newcommand{\OR}{\qquad\text{or}\qquad}
\newcommand{\chil}{\chi_{i,\La}}
\newcommand{\chjl}{\chi_{m+j,\La}}
\newcommand{\chkll}{\chi_{k,\Lambda^{(1)},\ldots,\Lambda^{(N)}}}
\newcommand{\chjll}{\chi_{m+j,\Lambda^{(1)},\ldots,\Lambda^{(N)}}}

\renewcommand{\theequation}{\arabic{section}.\arabic{equation}}
	
\newtheorem{thm}{Theorem}[section]
\newtheorem{lem}[thm]{Lemma}
\newtheorem{prop}[thm]{Proposition}
\newtheorem{cor}[thm]{Corollary}
\newtheorem{conj}[thm]{Conjecture}
\newtheorem*{mthm}{Main Theorem}
\newtheorem*{mthma}{Theorem A}
\newtheorem*{mthmb}{Theorem B}
	
\theoremstyle{definition}
\newtheorem{defin}[thm]{Definition}
	
\theoremstyle{remark}
\newtheorem{remark}[thm]{Remark}
\newtheorem{example}[thm]{Example}

\newtheorem{innercustomdef}{Definition}
\newenvironment{customdef}[1]
{\renewcommand\theinnercustomdef{#1}\innercustomdef}
{\endinnercustomdef}
	
\newcommand{\bth}{\begin{thm}}
\renewcommand{\eth}{\end{thm}}
\newcommand{\bpr}{\begin{prop}}
\newcommand{\epr}{\end{prop}}
\newcommand{\ble}{\begin{lem}}
\newcommand{\ele}{\end{lem}}
\newcommand{\bco}{\begin{cor}}
\newcommand{\eco}{\end{cor}}
\newcommand{\bde}{\begin{defin}}
\newcommand{\ede}{\end{defin}}
\newcommand{\bex}{\begin{example}}
\newcommand{\eex}{\end{example}}
\newcommand{\bre}{\begin{remark}}
\newcommand{\ere}{\end{remark}}
\newcommand{\bcj}{\begin{conj}}
\newcommand{\ecj}{\end{conj}}

\newcommand{\bal}{\begin{aligned}}
\newcommand{\eal}{\end{aligned}}
\newcommand{\beq}{\begin{equation}}
\newcommand{\eeq}{\end{equation}}
\newcommand{\ben}{\begin{equation*}}
\newcommand{\een}{\end{equation*}}
	
\newcommand{\bpf}{\begin{proof}}
\newcommand{\epf}{\end{proof}}

\newcommand{\Gaff}{\widehat{\mathfrak{g}}}   
\newcommand{\Gafft}{\widehat{\mathfrak{g}'}} 
\newcommand{\Glop}{L\mathfrak{g}}            
\newcommand{\Hlie}{\mathfrak{h}}             
\newcommand{\YB}{\mathcal{YB}}               
\newcommand{\CF}{\mathcal{F}}                
\newcommand{\CR}{\textbf{R}}                 
\newcommand{\BC}{\mathbb{C}}            
\newcommand{\BZ}{\mathbb{Z}}            
\newcommand{\Id}{\textrm{Id}}            
\newcommand{\CB}{\mathcal{B}}            
\newcommand{\Sm}{\mathbb{S}}             
\newcommand{\Rc}{\mathcal{R}}
\newcommand{\xiL}{\xi^{}_{\La}}
\newcommand{\xiLo}{\xi^{}_{\La^0}}
\newcommand{\etaL}{\eta^{}_{\La}}
\newcommand{\etaLo}{\eta^{}_{\La^0}}
\newcommand{\pa}[1]{\overline{#1}}	
	
\def\beql#1{\begin{equation}\label{#1}}
		
\title{Irreducibility of the tensor product of Yangian $\Y(\gl_{m|n})$ evaluation modules}
		
\author{Vyacheslav Futorny}
\address{Shenzhen International Center for Mathematics, Southern University of Science and Technology, Shenzhen, China}
\email{vfutorny@gmail.com}
\author{Zheng Li}
\address{School of Artificial Intelligence, Jianghan University, Wuhan, Hubei 430056, China}
\email{lz1994@jhun.edu.cn}
\author{Jian Zhang}
\address{School of Mathematics and Statistics,  Central China Normal University, Wuhan, Hubei 430079, China}
\email{jzhang@ccnu.edu.cn}

\thanks{{\scriptsize
\hskip -0.6 true cm MSC (2020): Primary: 17B37; Secondary: 81R10.
\newline Keywords: super Yangian, Gelfand-Tsetlin basis, evaluation module, tensor product.		
}}

\maketitle
		
\vspace{4 mm}		
\begin{abstract}

The evaluation homomorphism from the super Yangian $\Y(\gl_{m|n})$ to $\U(\gl_{m|n})$ induces a $\Y(\gl_{m|n})$-module structure on any finite dimensional simple $\U(\gl_{m|n})$-module $L(\la)$.
In this paper, we give  necessary and sufficient conditions for the tensor product of such evaluation
$\Y(\gl_{m|n})$-modules, $L_{g}(\la)\ot L_{h}(\ga)$, to be simple,  provided each of $\la$ and $\ga$ is either covariant tensor or essentially typical. Our proof is based on the existence of a Gelfand-Tsetlin basis for finite dimensional simple $\U(\gl_{m|n})$-modules with highest weights that belong to these two families: covariant tensor and essentially typical.

The obtained result is a super analogue of the Molev's result for the classical Yangian $\Y(\gl_n)$.  Combining this  with the binary property of tensor products of covariant evaluation modules, we obtain an irreducibility criterion for arbitrary tensor products of covariant evaluation modules.
\end{abstract}

\section{Introduction}

 Yangians is an important family of quantum groups deforming the universal enveloping algebras of current Lie algebras. They have deep connections with mathematical physics
 and numerous applications.

The Yangian $\Y(\gl_{n})$ of the general  linear algebra $\gl_{n}$ contains  a maximal commutative subalgebra \( A(\mathfrak{gl}_n) \) known as the  Gelfand-Tsetlin subalgebra.
It is generated by the centers of subalgebras in the chain
\[\Y(\mathfrak{gl}_1) \subset \Y(\mathfrak{gl}_2) \subset \cdots \subset \Y(\mathfrak{gl}_n)\]
and plays an important role in the representation theory of Yangians.
The generators of the center of each subalgebra $\Y(\mathfrak{gl}_k)$ are given by the coefficients of the corresponding quantum determinant \cite{MNO 1996,Mol 2007}.
A finite dimensional $\Y(\gl_{n})$-module is called \emph{tame} if the action of  the Gelfand-Tsetlin  subalgebra \( A(\mathfrak{gl}_n) \) is semisimple.
It was conjectured by Cherednik \cite{Che} and later proved by Nazarov and Tarasov \cite{NT 1998} that a simple $\Y(\gl_{n})$-module is tame if and only if it
is obtained by pulling back through some automorphism $\om_f$ (cf. Section 3, \eqref{omf}) from the tensor product
\begin{equation}\label{tensor product}
V_{h_1}(\la_1\slash\mu_1)\otimes V_{h_2}(\la_2\slash\mu_2)\otimes \cdots\otimes V_{h_k}(\la_k\slash\mu_k)
\end{equation}
of \emph{elementary} (or \emph{skew}) \( \Y(\mathfrak{gl}_n) \)-modules, for some skew Young diagrams \( \la_1\slash\mu_1, \ldots, \la_k\slash\mu_k \) and some complex numbers \( h_1, \ldots, h_k \) such that \( h_i - h_j \notin \mathbb{Z} \) when \( i \neq j \).
 In particular, under the latter restriction
 the tensor product  \eqref{tensor product} of the elementary \( \Y(\mathfrak{gl}_n) \)-modules is always simple.

Nazarov and Tarasov also proved that  a general tensor product \begin{equation}
V_{z_1}(\la_1\slash\mu_1)\otimes V_{z_2}(\la_2\slash\mu_2)\otimes \cdots\otimes V_{z_k}(\la_k\slash\mu_k)
\end{equation} is simple if and only if
\[V_{z_i}(\la_i\slash\mu_i)\otimes V_{z_j}(\la_j\slash\mu_j)\] is simple for any increasing pair \((i, j)\), \(i, j=1, \ldots, k\) .
On the other hand,  Molev \cite{Mol 2002} showed that the  tensor product
of  two highest weight $\gl_{n}$-modules $L(\lambda)$ and $L(\mu)$  is a simple \(\Y(\mathfrak{gl}_n)\)-module if and only if the highest weights \(\lambda\) and \(\mu\) satisfy the non-crossing condition. Together with the binary condition of Nazarov and Tarasov it gives the irreducibility of tensor product of evaluation modules.

The super  analogue of the Yangian $\Y(\gl_n)$, the super Yangian $\Y(\gl_{m|n})$, was introduced by Nazarov \cite{Naz 1991}.
Finite dimensional irreducible representations of $\Y(\gl_{m|n})$ were classified by Zhang in \cite{Zhang 1995, Zhang 1996}.
The evaluation map from $\Y(\gl_{m|n})$ to $\U(\gl_{m|n})$  allows one to view $\gl_{m|n}$-modules as $\Y(\gl_{m|n})$-modules, giving a family of evaluation modules.
Any finite dimensional irreducible representation of $\Y(\gl_{m|n})$ is a subquotient of tensor products of evaluation modules.

In contrast to the case of $\Y(\gl_n)$, not every finite dimensional simple module $L(\lambda)$ of $\gl_{m|n}$ is tame.  There are two known types of tame modules: \emph{covariant tensor} modules and \emph{essentially typical} modules.
In \cite{FLZ 2026}, the authors studied simple quotients of the submodules of  tensor products of covariant evaluation modules generated by a tensor product
of highest weight vectors,  and gave necessary and sufficient conditions for such modules to be tame. This generalized the earlier work of Nazarov and Tarasov \cite{NT 1998} for
$\Y(\mathfrak{gl}_m)$ to the super case.

Provided $\la$ is covariant, $L(\la)$ can also be obtained as a submodule of a tensor product of natural representations. This makes it possible to adapt the strategy of Nazarov and Tarasov to the super Yangian setting, as observed in \cite{LM 2021}.
Despite the close analogy between Yangians and super Yangians, no irreducibility criterion for tensor products of evaluation modules over $\Y(\gl_{m|n})$ was known so far. The main difficulty lies in the fact that the techniques used in the non-super case do not carry over directly to the super setting. For example, the proof of the sufficiency of Molev's criterion in the non-super case relies heavily on quantum minors. But when it comes to the super case, the quantum minors should be replaced by the Berezinian minors. However, the coproduct formulas for Berezinian minors are considerably more involved. Furthermore, unlike quantum minors, Berezinian minors do not interact well with the weight structure. To overcome these difficulties, we develop a different approach based on a refined analysis of Gelfand--Tsetlin bases.

As a consequence, we show that the Molev's non-crossing condition admits a natural extension to the super Yangian setting. Namely, we have the following main result.

\begingroup
\renewcommand{\thethm}{A}
\bth[Proposition \ref{prop:shif}, Theorem \ref{thm:scon}, Theorem \ref{thm:ncon}]
Let each of $\la$ and $\ga$ be either essentially typical or covariant $\gl_{m|n}$-weight and $g,h\in\BC$. Then the $\Y(\gl_{m|n})$-module $L_g(\la)\ot L_h(\ga)$ is simple if and only if  $\la+g$ and $\ga+h$ are super non-crossing.
\eth
\endgroup

Combined with the binary property, this theorem yields an irreducibility criterion for arbitrary finite tensor products of covariant evaluation modules.

\begingroup
\renewcommand{\thethm}{B}
\bth[Corollary \ref{coro:multi}]
Let $\la_i$ be a covariant weight and  $h_i\in \BC$ for $1\leq i\leq k$.
Then the $\Y(\gl_{m|n})$-module
\[
L_{h_1}(\la^{(1)})\ot L_{h_2}(\la^{(2)})\ot\cdots\ot L_{h_k}(\la^{(k)}),
\]
is simple if and only if  $\la_i+h_i$ and $\la_j+h_j$ are super non-crossing for $1\leq i< j\leq k$.
\eth
\endgroup

Furthermore,  the irreducibility criterion established in the present paper shows that, in many cases, a tensor product of covariant evaluation modules is simple, and hence isomorphic to a submodule generated by the tensor product of highest weight vectors. Consequently, this makes the construction in \cite{FLZ 2026} more transparent.

The paper is organized as follows.
In Section 2, we recall some preliminaries of Lie superalgebra $\gl_{m|n}$, especially the Gelfand-Tsetlin basis for covariant and essentially typical $\gl_{m|n}$-modules. In Section 3, we recall the definition of super Yangian $\Y(\gl_{m|n})$ following \cite{Gow 2005, Gow 2007}. Moreover some facts about the representations of $\Y(\gl_{m|n})$ are presented for later use.
In Section 4, we first prove using the Gelfand-Tsetlin basis that the super non-crossing condition is a sufficient condition for $L_g(\la)\ot L_h(\ga)$ to be simple.
After that we show the super non-crossing condition is also necessary  for the irreducibility of
$L_g(\la)\ot L_h(\ga)$.

\section{Lie superalgebra $\gl_{m|n}$}
Throughout the paper, we work over $\BC$. A super vector space $V=\BC^{m|n}$ is a ${\BZ}_2$-graded vector space.
Vectors in the $\bar{0}$-graded part $({\BC}^{m|n})_{\bar{0}}={\BC}^m=\sum_{i=1}^m {\BC}
v_i$ are called even, while those in the $\bar{1}$-graded part $({\BC}^{m|n})_{\bar{1}}={\BC}^n=\sum_{j=m+1}^{m+n}{\BC} v_j$ are called odd.
We define the parity $\bi$ of $i$ by
\beql{parity}
\bi =
\begin{cases}
0, & \text{if } i \leq m, \\
1, & \text{if } i > m.
\end{cases}
\eeq

The Lie superalgebra $\mathfrak{gl}_{m|n}$ has a basis $\{E_{ij}\mid 1\le i,j\le m+n\}$, where the parity of $E_{ij}$ is $\bar i+\bar j$. The supercommutator is given by
\[
[E_{ij},E_{kl}] = \delta_{jk}E_{il} - (-1)^{(\bar i+\bar j)(\bar k+\bar l)}\delta_{il}E_{kj}.
\]

Denote by $\U(\gl_{m|n})$ the universal enveloping superalgebra of $\gl_{m|n}$.
The Cartan subalgebra $\h$ of $\gl_{m|n}$ is spanned by $E_{ii}$, $1\leq i\leq m+n$. Let $\ve_i$ ($1\leq i\leq m+n$) be a basis of $\h^*$ (the dual space of $\h$) such that $\ve_i(E_{jj})=\delta_{ij}$.
There is a bilinear form $(\ ,\ )$ on $\h^*$ given by $(\ve_i,\ve_j)=(-1)^{\bi}\delta_{ij}$.
Define the simple roots $\alpha_i:=\ve_i-\ve_{i+1}$, for $1\leq i\leq m+n-1$.

Given a tuple $\la$ of complex numbers of the form $\la=(\la_{1},\la_{2},\ \ldots,\la_{m+n})$, the simple highest weight module $L(\lambda)$ of the Lie superalgebra $\gl_{m|n}$ with the highest weight $\lambda$ is generated by a nonzero vector $\xi$ satisfying the conditions:

\beq
\begin{split}
&E_{ii} \xi=\lambda_i \xi, 1\leq i\leq m+n, \\
&E_{ij} \xi=0,  1\leq i<j\leq m+n-1.
\end{split}
\eeq

Recall that the  weight $\lambda$ is dominant if
\beq
\la_i-\la_{i+1}\in \BZ_{\geq 0}, \text{ for }i=1, \ldots,m+n-1,  \quad i\neq m.
\eeq

\subsection{Gelfand-Tsetlin basis (essentially typical case)}\label{sec:GT basis}
The dominant weight $\la$ (and the module $L(\la)$) is \emph{typical} if
$(\lambda+\rho,\varepsilon_i-\varepsilon_j)\neq 0$ for  $1\leq i\leq m,1\leq j\leq n$, where $2\rho$ is the sum of all positive roots of $\gl_{m|n}$.

Set
\beq
l_{i}=\la_{i}-i+1, (1\leq i\leq m); \quad l_{j}=-\la_{j}+j-2m, (m+1\leq j\leq m+n).
\eeq
Since $(\lambda+\rho,\ve_i-\ve_j)=l_{i}-l_{j}$,
the condition for typicality is equivalent to $l_{i}\neq l_{j}$ for  $1\leq i\leq m <j\leq m+n. $

For $a, b \in \CC$ such that $b-a$ is nonnegative integer, we set
\[
[a;b]=\{a+s \ | 0\leq s\leq b-a\}.
\]

The dominant weight $\la$ (and the module $L(\la)$) is called \emph{essentially typical} if
\beq
\{l_1,l_2,\ldots,l_m\}\cap [l_{m+1}; l_{m+n}]=\emptyset.
\eeq
Clearly,  every essentially typical weight is typical.

Our main combinatorial device is an array of complex numbers
$\Lambda=(\lambda_{ij})$ presented in the following form:

\beql{GT pattern}
\begin{array}{cccccccc}
\la_{m+n,1}   &  \cdots      & \la_{m+n,m} & \la_{m+n,m+1} & \cdots & \la_{m+n-1,m+n} & \la_{m+n,m+n} \\
\la_{m+n-1,1} &  \cdots     & \la_{m+n-1,m} & \la_{m+n-1,m+1} & \cdots & \la_{m+n-1,m+n-1} & \\
\vdots    &  \vdots & \vdots & \vdots   & \reflectbox{$\ddots$}    \\
\la_{m+1,1} &  \cdots  & \la_{m+1,m} & \la_{m+1,m+1}    \\
\la_{m1} &  \cdots  & \la_{mm}                           \\
\la_{m-1,1} &  \cdots                             \\
\vdots  &  \reflectbox{$\ddots$}                              \\
\la_{11}\\
\end{array}
\eeq
We will call  $\Lambda$ a  \emph{Gelfand-Tsetlin pattern} or simply a pattern.  Given $\Lambda=(\lambda_{ij})$
we set

\beq
l_{ki}=\la_{ki}-i+1, (1\leq i\leq m); \quad l_{kj}=-\la_{kj}+j-2m, (m+1\leq j\leq k).
\eeq

\bth[\cite{Palev 1989a}, \cite{Palev 1989b}]\label{thm:GT basis}
Let $\la$ be an essentially typical highest weight. Then $L(\la)$
admits a basis $\xi_{\La}$ parameterized by all patterns $\La$ satisfying
following conditions:
\begin{enumerate}
\item $\lambda_{m+n,i}=\lambda_i, \quad 1\leq i\leq m+n$;
\item  $\la_{k,i}-\la_{k-1,i}\equiv\theta_{k-1,i}\in\{0,1\}, 1\leq i\leq m;m+1\leq k\leq m+n$;
\item $\la_{ki}-\la_{k,i+1}\in \mathrm{Z}_{\geq 0}, 1\leq i\leq m-1;m+1\leq k\leq m+n-1$;
\item  $\la_{k+1,i}-\la_{ki}\in \mathrm{Z}_{\geq 0}$ {\it and} $\la_{k,i}-\la_{k+1,i+1}\in \mathrm{Z}_{\geq 0},$
$1\leq i\leq k\leq m-1$ \text{ or } $m+1\leq i\leq k\leq m+n-1$.
\end{enumerate}
The action of the generators of $\gl_{m,n}$
is given by the formulas
\beql{action h}
E_{kk}\xi_{\La}=\left(\sum_{i=1}^{k}\la_{kj}-\sum_{j=1}^{k-1}\la_{k-1,j}\right)\xi_{\La},\quad 1\leq k\leq m+n;
\\
\eeq
	
\beq
\noindent E_{k,k+1}\xi_{\La}=-\sum_{i=1}^{k}\frac{\Pi_{j=1}^{k+1}(l_{k+1,j}-l_{ki}) }
{\Pi_{j\neq i,j=1}^{k} (l_{kj}-l_{ki}) }\xi_{\La+\delta_{ki}},
\quad1\leq k\leq m-1;
\eeq
	
\beq
E_{k+1,k}\xi_{\La}=\sum_{i=1}^{k}\frac{\Pi_{j=1}^{k-1}(l_{k-1,j}-l_{ki})}
{\Pi_{j\neq i,j=1}^{k}(l_{kj}-l_{ki})}\xi_{\La-\delta_{ki}},
\quad  1\leq k\leq m-1;
\eeq
	
\beq
\begin{split}
E_{m,m+1}\xi_{\La}=\sum_{i=1}^{m}\theta_{mi}(-1)^{i-1}(-1)^{\theta_{m1}+\ldots+\theta_{m,i-1}}
\\
\times  \frac{\Pi_{1\leq j< i} (l_{mj}-l_{mi}-1)}
{\Pi_{i<j\leq m} (l_{mj}-l_{mi})
\Pi_{j\neq i,j=1}^{m}(l_{m+1,j}-l_{mi}-1)}
\xi_{\La+\delta_{mi}},
\\
\end{split}
\eeq

\beq
\begin{split}
E_{m+1,m} \xi_{\La}=\sum_{i=1}^{m}(1-\theta_{mi})(-1)^{i-1}(-1)^{\theta_{m1}+\ldots+\theta_{m,i-1}}
\\
\times  \frac{(l_{m,i}-l_{m+1,m+1})\Pi_{ i<j\leq m} (l_{mj}-l_{mi}+1)\Pi_{j=1}^{m-1}(l_{m-1,j}-l_{mi})}
{\Pi_{1\leq j< i} (l_{mj}-l_{mi})} \xi_{\La-\delta_{mi}},
\end{split}
\eeq

\beq
\begin{split}
E_{k,k+1} \xi_{\La}=\sum_{i=1}^{m}\theta_{ki}(-1)^{\theta_{k1}+\ldots+\theta_{k,i-1}
+\theta_{k-1,i+1}+\ldots+\theta_{k-1,m}}(1-\theta_{k-1,i})
\\
\times
\prod_{j\neq i,j =1}^{m}\left(\frac{l_{kj}-l_{ki}-1}{l_{k+1,j}-l_{ki}-1}\right)
\xi_{\La+\delta_{ki}}
\\
\\
-\sum_{i=m+1}^{k}
\Pi_{j=1}^{m}\left(\frac{(l_{kj}-l_{ki})(l_{kj}-l_{ki}+1)}{(l_{k+1,j}-l_{ki})(l_{k-1,j}-l_{ki}+1)} \right)
\\
\times
\frac{\Pi_{j=m+1}^{k+1}(l_{k+1,j}-l_{ki})}
{\Pi_{j\neq i,j=m+1}^{k} (l_{kj}-l_{ki})}\xi_{\La+\delta_{ki}} ,
\quad  m+1\leq k\leq m+n-1;
\\
\end{split}
\eeq
	
\beql{action f}
\begin{split}
E_{k+1,k}\xi_{\La}=
\sum_{i=1}^{m}\theta_{k-1,i}(-1)^{\theta_{k1}+\ldots+\theta_{k,i-1}+\theta_{k-1,i+1}+\ldots+\theta_{k-1,m}}(1-\theta_{ki})
\\
\times
\prod_{j\neq i=1}^{m}\left(\frac{l_{kj}-l_{ki}+1}{l_{k-1,j}-l_{ki}+1}\right)
\frac{\Pi_{j=m+1}^{k+1}(l_{k+1,j}-l_{ki})\Pi_{j=m+1}^{k-1}(l_{k-1,j}-l_{ki}+1)}
{\Pi_{j=m+1}^{k} (l_{kj}-l_{ki})(l_{kj}-l_{ki}+1)}\xi_{\La-\delta_{ki}}
\\
+ \sum_{i=m+1}^{k}
\frac{\prod_{j=m+1}^{k-1}(l_{k-1,j}-l_{ki})}{\prod_{j\neq i,j =m+1}^{k}(l_{k,j}-l_{ki})}
\xi_{\La-\delta_{ki}} \quad
m+1\leq k\leq m+n-1.
\end{split}
\eeq
The arrays $\La\pm \delta_{ki}$ are obtained from $\La$ by replacing $\la_{ki}$ by $\la_{ki}\pm1$. We assume  that
$\xi_{\La}=0$ if the array $\La$ does not satisfy the conditions of the theorem.
\eth

\subsection{Gelfand-Tsetlin basis (covariant case)}
The Lie superalgebra $\gl_{m|n}$ has a natural representation on $\BC^{m|n}$ called the vector representation, such that $E_{ij}v_k=\delta_{jk}v_i$. The highest weight of $\BC^{m|n}$ is the tuple $(1,0,\dots,0)$.

We call the weight $\la$ (and the module $L(\la)$) \emph{covariant},
if $\la$ satisfies:
\begin{enumerate}
\item $\la_1,\ldots,\la_{m+n}$ are all nonnegative integers,
\item $\la_1\geq \ldots\geq \la_m$ and $\la_{m+1}\geq\ldots\geq\la_{m+n}$,
\item $\sharp\{\la_{m+j}>0|1\leq j\leq n\}\leq \la_m$.
\end{enumerate}
Note that in this case $L(\la)$ is a submodule of $(\BC^{m|n})^{\otimes |\la|}$, where $|\la|=\sum_{i=1}^{m+n}\la_i$.

\bth[\cite{Mol 2011}, \cite{SV 2010}]\label{thm:GT basis2}
Let $\la$ be a covariant weight. Then $L(\la)$
admits a basis $\xi_{\La}$ parameterized by all patterns $\La$ satisfying the
following conditions:
\begin{enumerate}
 \item $\lambda_{m+n,i}=\lambda_i, \quad 1\leq i\leq m+n$;
 \item  $\la_{k,i}-\la_{k-1,i}\equiv\theta_{k-1,i}\in\{0,1\}, 1\leq i\leq m;m+1\leq k\leq m+n$;
 \item $\la_{ki}-\la_{k,i+1}\in \mathrm{Z}_{\geq 0}, 1\leq i\leq m-1;m+1\leq k\leq m+n-1$;
 \item  $\la_{k+1,i}-\la_{ki}\in \mathrm{Z}_{\geq 0}$ {\it and} $\la_{k,i}-\la_{k+1,i+1}\in \mathrm{Z}_{\geq 0},$
$1\leq i\leq k\leq m-1$ \text{ or } $m+1\leq i\leq k\leq m+n-1$;
  \item
$m+1\le k\le m+n$:
$\lambda_{km}\ge \#\{i:\lambda_{ki}>0,\; m+1\le i\le k\};$
  \item
if $\lambda_{ m+1,m}=0$, then $\theta_{mm}=0$.
\end{enumerate}
The action of generators of $\gl_{m|n}$ on this basis satisfies the same formulas given in Theorem \ref{thm:GT basis}.
\eth

\bre
Essentially typical representations of the Lie superalgebra $\gl(m|n)$ with a Gelfand-Tsetlin basis were constructed in \cite{Palev 1989a}, \cite{Palev 1989b}. This construction was later generalized to \emph{covariant tensor} modules in \cite{SV 2010}.
Another explicit construction of covariant tensor modules, based on super Young tableaux, was given by Molev \cite{Mol 2011}.
In the following we will use the formulas from \cite{FSZ 2021},
which are obtained by a suitable modification of the formulas of Palev \cite{Palev 1989b} and Stoilova and Van der Jeugt \cite{SV 2010}.
\ere

If a $\gl_{m|n}$-weight $\la=(\la_1,\ldots,\la_{m+n})$ satisfies the condition that $\la+h=(\la_1+h,\ldots,\la_m+h,\la_{m+1}-h,\ldots,\la_{m+n}-h)$ is covariant for some $h\in\BC$, then we call $\la$ a \emph{shifted covariant weight}.
It follows immediately that for such $\la$, the simple highest weight module $L(\la)$ admits a Gelfand-Tsetlin basis parametrized by patterns $\La=(\la_{ij})$ such that $\La+h=(\la_{ij}+(-1)^{\bj}h)$ satisfies the conditions of Theorem \ref{thm:GT basis2} for $L(\la+h)$. At the same time, it is clear that a shift of essentially typical weight is still essentially typical.

We call a pattern \emph{admissable} associated with a essentially typical or shifted covariant $\la$ if (after a shift) it satisfies the conditions listed in Theorem \ref{thm:GT basis} or Theorem \ref{thm:GT basis2}, respectively.
The set of all admissable patterns associated with $\la$ is denoted by $\mathscr{S}_{\la}$.

\section{Super Yangian $\Y(\gl_{m|n})$}
The Yangian $Y (\mathfrak{gl}_{m|n})$
is the $\mathbb{Z}_2$-graded
associative algebra over $\mathbb{C}$ with generators
\[
\{t_{ij}^{(r)}\,| \; 1\le i,j \le m+n; r\ge 1\}
\]
and defining relations
\beql{comr}
[t_{ij}^{(r)}, t_{kl}^{(s)}] = (-1)^{\pa{i}\,\pa{j} + \pa{i}\,\pa{k} + \pa{j}\,\pa{k}}
\sum_{p=0}^{\mathrm{min}(r,s) -1}(t_{kj}^{(p)} t_{il}^{(r+s-1-p)} - t_{kj}^{(r+s-1-p)}t_{il}^{(p)}).
\eeq
where $\pa{i}$ is the parity of the index $i$.
We take $\pa{i} = 0$ for $i\le m$; and $\pa{i}=1$ for $i \ge m+1$.
(We write square brackets for the super-commutator).
We define the formal power series
\beq
t_{ij} (u) = \delta_{ij} + t_{ij}^{(1)} u^{-1} + t_{ij}^{(2)}u^{-2} + \ldots\in\Y(\gl_{m|n})[[u^{-1}]],
\eeq
then
\beql{comr2}
[t_{ij}(u),\;t_{kl}(v)] \;=\; \frac{(-1)^{\pa{i}\,\pa{j}+\pa{i}\,\pa{k}+\pa{j}\,\pa{k}}}{(u-v)}(t_{kj}(u)
t_{il} (v) - t_{kj} (v) t_{il} (u) ).
\eeq
Define the generating matrix as
\beql{T}
T(u) = \sum_{i,j =1}^{m+n} E_{ij} \otimes  t_{ij} (u) (-1)^{\pa{j} (\pa{i} +1)}\in \End(\BC^{m|n})\ot\Y(\gl_{m|n})[[u^{-1}]]
\eeq
where $E_{ij}$ is the standard elementary matrix.
Consider the permutation matrix
\[
P = \sum_{i,j =1}^{m+n} E_{ij} \otimes E_{ji} (-1)^{\pa{j}}.
\]
The rational function
\[
R(u) = 1 - \frac{P}{u}
\]
satisfies the Yang-Baxter equation
\[
R_{12}(u)R_{13}(u+v)R_{23}(v)=R_{23}(v)R_{13}(u+v)R_{12}(u).
\]
Then,
the defining relations may be expressed by the matrix product
\beql{RTT}
R(u-v) T_1 (u) T_2 (v) = T_2 (v) T_1 (u) R(u-v),
\eeq
which holds in $\End(\BC^{m|n})^{\ot 2}\ot\Y(\gl_{m|n})[[u^{-1}]]$.

The Yangian $\Y(\gl_{m|n})$ is a Hopf algebra with coproduct
\beql{copro1}
\Delta : t_{ij}(u)  \mapsto \sum_{k=1}^{m+n} t_{ik}(u)  \otimes t_{kj}(u),
\eeq
antipode $S: T(u) \mapsto T(u)^{-1}$ and counit $\epsilon : T(u) \mapsto 1$.

Another coproduct $\bar{\Delta}$ on $\Y(\gl_{m|n})$ is defined by
\beql{copro2}
\bar{\De}:t_{ij}(u)  \mapsto \sum_{k=1}^{m+n} (-1)^{(\bi+\bk)(\bj+\bk)}t_{kj}(u)  \otimes t_{ik}(u).
\eeq
 Using \eqref{RTT}, it is easy to see that
the map $\om_f:\Y(\gl_{m|n})\rightarrow \Y(\gl_{m|n})$ defined by
\beql{omf}
\om_f:t_{ij}(u)=f(u)t_{ij}(u)
\eeq
for $1\leq i,j\leq m+n$ and $f(u)\in \BC[[u^{-1}]]$ is an automorphism.

\bpr[\cite{Naz 2020}]\label{pro:anti-auto}
The map $\st:\Y(\gl_{m|n})\rightarrow \Y(\gl_{m|n})$ defined by
\beql{st}
\st(t_{ij}(u))=(-1)^{\bi(\bj+1)}t_{ji}(u)
\eeq
for $1\leq i,j\leq m+n$, is an anti-automorphsim.
\epr

\subsection{Berezinian}
Throughout this paper, we use the following notation for  the inverse of $T(u)$
\[
T(u)^{-1}=\left(t'_{ij}(u)\right)_{i,j=1}^{m+n}.
\]
Denote by $\mathcal{S}_k$ the $k$-th symmetric group.
\bde
The Berezinian of $\Y(\gl_{m|n})$ is defined by
\begin{equation}\label{bere}
\begin{split}
B(u)=&\sum_{\si\in \mathcal{S}_m}\sgn\si t_{\si(1)1}(u)t_{\si(2)2}(u-1)\cdots t_{\si(m)m}(u-m+1)\\
&\times \sum_{\tau\in \mathcal{S}_n}\sgn\tau t'_{m+1,m+\tau(1)}(u-m+1)t'_{m+2,m+\tau(2)}(u-m+2)\cdots t'_{m+n,m+\tau(n)}(u-m+n)
\end{split}
\end{equation}
\ede

The following result was conjectured by Nazarov\cite{Naz 1991} and later proved by Gow\cite{Gow 2007}.
\bth[{\cite{Gow 2007}}]\label{thm:B}
The coefficients of $B(u)$ generate the center of $\Y(\gl_{m|n})$.
\eth
For the $(m+n)\times (m+n)$ matrix $T(u)$, let $T^{(k)}(u)$ be the submatrix consisting of its first $k$ rows and first $k$ columns.
For $1\leq i\leq m$ and $1\leq j\leq n$, we define Berezinian minors
\[
B_i(u)=\sum_{\si\in \mathcal{S}_i}\sgn\si t_{\si(1)1}(u)t_{\si(2)2}(u-1)\cdots t_{\si(i)i}(u-i+1),
\]
\[
\begin{split}
B_{m+j}(u)=&\sum_{\si\in \mathcal{S}_m}\sgn\si t_{\si(1)1}(u)t_{\si(2)2}(u-1)\cdots t_{\si(m)m}(u-m+1)\\
&\times \sum_{\tau\in \mathcal{S}_j}\sgn\tau (T^{(m+j)}(u-m+1)^{-1})_{m+1,m+\tau(1)}\cdots (T^{(m+j)}(u-m+j)^{-1})_{m+j,m+\tau(j)}
\end{split}
\]
The coefficients of $B_i(u)$ and $B_{m+j}(u)$ generate the centers of $\Y(\gl_i)$ and $\Y(\gl_{m|j})$ respectively.

Considering the symmetrizer $\sum_{\si\in\mathcal{S}_k}\si$ and anti-symmetrizer $\sum_{\si\in\mathcal{S}_k}\sgn(\si)\si$ in $\mathcal{S}_k$,
their image under the natural action of $\mathcal{S}_k$ on $(\mathbb{C}^{m|n})^{\ot k}$ are denoted
by $A_k$ and $H_k$, respectively.
Then the Berezinian $B(u)$ equals the coefficient of $E_{11}\ot E_{22}\ot\cdots\ot E_{m+n,m+n}$ in
\[
A_mT_1(u)\cdots T_m(u-m+1)H_{n}^{\circ}T_{m+1}(u-m+1)^*\cdots T_{m+n}(u-m+n)^*,
\]
where $T(u)^*=(T(u)^{-1})^{\text{st}}$ is the supertranspose of $(T(u))^{-1}$ and $H_{n}^{\circ}$ is the symmetrizer in $\mathcal{S}_n$ acting
on the
tensor factors of  $(\mathbb{C}^{m|n})^{\otimes m+n}$ indexed  by $\{m+1,\ldots, m+n\}$\cite{Naz 1991,Naz 2020}.

Similarly, $B_i(u)$ equals the coefficient of $E_{11}\ot E_{22}\ot\cdots\ot E_{ii}$ in
\[
A_iT_1(u)\cdots T_i(u-i+1),
\]
and
$B_{m+j}(u)$ equals the coefficient of $E_{11}\ot E_{22}\ot\cdots\ot E_{m+j,m+j}$ in
\[
A_mT_1(u)\cdots T_m(u-m+1)H_{j}^{\circ}(T^{(m+j)}(u-m+1)^*)_{m+1}\cdots (T^{(m+j)}(u-m+j)^*)_{m+j},
\]
where $T^{(k)}(u)^*=(T^{(k)}(u)^{-1})^{\text{st}}$ and $H_{j}^{\circ}$ is the symmetrizer $H_j$ acting
on the tensor factors  of  $(\mathbb{C}^{m|n})^{\otimes m+j}$ indexed by $\{m+1,\ldots, m+j\}$.

Let us consider the space $\End(\mathbb{C}^{m|n})^{\ot k}\ot \Y(\gl_{m|n})[[u^{-1}]]\ot \Y(\gl_{m|n})[[u^{-1}]]$.
For a matrix $A=\sum e_{ij}\ot a_{ij}$, let
\[
A_{p[q]}=\sum 1^{\ot (p-1)}\ot E_{ij}\ot 1^{\ot (k-p)}\ot 1^{\ot (q-1)}\ot a_{ij}\ot 1^{\ot (2-q)},
\]
where $1\leq p\leq k$ and $1\leq q\leq 2$.
Then we can rewrite \eqref{copro1} as
\[
\Delta(T(u))=T(u)_{1[1]}T(u)_{1[2]}.
\]
It also implies that
\[
\Delta(T(u)^{-1})=(T(u)^{-1})_{1[2]}(T(u)^{-1})_{1[1]}.
\]

Let $A$ be any square matrix of size $m+n\times m+n$.
For any subsets $I=\{i_1,\ldots,i_p\}$, $J=\{j_1,\ldots,j_q\}$ of $\{1,2,\ldots,m+n\}$,
we denote by $A^I_J$ the matrix whose $ab$-th entry is $A_{i_a j_b}$.
Then
\[
\Delta(T^{(m+j)}(u))=(T(u)^{\mathbb{M}+\mathbb{J}}_{\mathbb{M}+\mathbb{N}})_{1[1]} (T(u)^{\mathbb{M}+\mathbb{N}}_{\mathbb{M}+\mathbb{J}})_{1[2]}
\]
where $\mathbb{M}+\mathbb{J}=\{1,2,\ldots,m+j\}$ and $\mathbb{M}+\mathbb{N}=\{1,2,\ldots,m+n\}$.
Similarly, there is
\[
\Delta(T^{(m+j)}(u)^*)=((T(u)^*)^{\mathbb{M}+\mathbb{J}}_{\mathbb{M}+\mathbb{N}})_{1[1]}((T(u)^*)^{\mathbb{M}+\mathbb{N}}_{\mathbb{M}+\mathbb{J}})_{1[2]}.
\]

As a result, $\Delta(B_i(u))$ is the coefficient of $E_{11}\ot E_{22}\ot\cdots\ot E_{ii}$ in
\beql{com1}
A_iT_{1[1]}(u)\cdots T_{i[1]}(u-i+1)A_iT_{1[2]}(u)\cdots T_{i[2]}(u-i+1),
\eeq
and
$\Delta(B_{m+j}(u))$ is the coefficient of $E_{11}\ot E_{22}\ot\cdots\ot E_{m+j,m+j}$ in
\beql{com2}
\begin{split}
&A_mT_{1[1]}(u)\cdots T_{m[1]}(u-m+1)H_{j}^{\circ}((T(u)^*)^{\mathbb{M}+\mathbb{J}}_{\mathbb{M}+\mathbb{N}})_{m+1[1]}\cdots ((T(u)^*)^{\mathbb{M}+\mathbb{J}}_{\mathbb{M}+\mathbb{N}})_{m+j[1]}\\
\times &A_mT_{1[2]}(u)\cdots T_{m[2]}(u-m+1)H_{j}^{\circ}((T(u)^*)^{\mathbb{M}+\mathbb{N}}_{\mathbb{M}+\mathbb{J}})_{m+1[2]}\cdots ((T(u)^*)^{\mathbb{M}+\mathbb{N}}_{\mathbb{M}+\mathbb{J}})_{m+j[2]},
\end{split}
\eeq

\subsection{Evaluation module}
There is an injective homomorphism $\iota:U(\gl_{m|n})\rightarrow \Y(\gl_{m|n})$ given by
\[
\iota: E_{ij}\rightarrow t_{ij}^{(1)}(-1)^{\bi},
\]
and a surjective homomorphism $\pi_{m|n}:\Y(\gl_{m|n})\rightarrow U(\gl_{m|n})$ defined as follows:
\[
\pi_{m|n}: t_{ij}(u)\rightarrow \delta_{ij}+E_{ij}(-1)^{\bi}u^{-1},
\]
which is called the evaluation homomorphism.
Under the homomorphism $\iota$, we may regard $E_{ij}$ as an element of $\Y(\gl_{m|n})$.


The assignment
\[
t_{ij}(u)\mapsto t_{ij}(u+h),
\]
where $h\in \BC$, defines an automorphism of $\Y(\gl_{m|n})$, called a shift automorphism.
The composition of the evaluation homomorphism and a shift automorphism allows one to regard a $\gl_{m|n}$-module as a $\Y(\gl_{m|n})$-module.
Such module is called an \emph{evaluation module} of $\Y(\gl_{m|n})$. As a result, for any $h\in \BC$ the simple highest weight module $L(\la)$ of $\gl_{m|n}$ becomes a $\Y(\gl_{m|n})$-module denoted by $L_h(\la)$. When $h=0$, we write $L_h(\la)$ simply as $L(\la)$.

For any two $\gl_{m|n}$-weights $\la$ and $\ga$, $L(\la)\ot L(\ga)$ can be viewed as a $\Y(\gl_{m|n})$-module through the coproduct $\De$ given in \eqref{copro1}. Moreover, it can also be regarded as a $\Y(\gl_{m|n})$-module through $\bar{\De}$ given in \eqref{copro2}. In the latter case, we denote it by
$\overline{L(\la)\ot L(\ga)}$.

For a $\Y(\gl_{m|n})$-module $M$, its dual space $M^*$ is a $\Y(\gl_{m|n})$-module through the anti-automorphism $\st$ given in Proposition \ref{pro:anti-auto}.
If $N$ is a proper submodule of $M$, it is clear that
\[
\text{Ann}{N}=\{f\in M^*|f(\eta)=0 \text{ for all } \eta\in N\}
\]
is a proper submodule of $M^*$.
Conversely, if $L$ is a proper submodule of $M^*$, then
\[
\text{Ker}{L}=\{\eta\in M|f(\eta)=0 \text{ for all }f\in L\}
\]
is a  proper submodule of $M^*$.
It implies that $M$ is simple if and only if $M^*$ is simple.

 The following proposition holds.
\bpr\label{pro:dual} For any  $\gl_{m|n}$-weights $\la$ and $\ga$
we have a $\Y(\gl_{m|n})$-module  isomorphism:
\[(L(\la)\ot L(\ga))^*\simeq
\overline{L(\la)\ot L(\ga)}.
\]
\epr

\bpf
We first show that $L(\la)^*\cong L(\la)$.
Due to the above observation, $L(\la)^*$ is simple.
Let $\xi$ be a highest weight vector of $L(\la)$. Then
 the vector $f_{\xi}$ defined by
\[
f_{\xi}(\xi)=1,\quad f_{\xi}(\eta)=0 \text{ for all }\eta\in L(\la)_{\mu}\text{ with }\mu\neq\la
\]
is a highest weight vector of $L(\la)^*$.
 Finding the weight of $f_{\xi}$ we obtain the claim.

Since $(L(\la)\ot L(\ga))^*$ can be naturally identified with $L(\la)^*\ot L(\ga)^*$ as a vector space, the result
follows from the following equality
\[
\De\circ\st=(\st\ot\st)\circ\bar{\De}.
\]
\epf

In this paper, we consider the irreducibility of the $\Y(\gl_{m|n})$-module $L_g(\la)\ot L_h(\ga)$ when each of $\la$ and $\ga$ is either essentially typical or covariant.

\subsection{Drinfeld polynomials}
Let $L(\la(u))$ be the simple highest weight module of $\Y(\gl_{m|n})$, i.e. it is generated by a highest weight vector satisfying that
\begin{enumerate}
\item $t_{ij}(u)v=0$ for $1\leq i<j\leq m+n$,
\item $t_{ii}(u)v=\la_i(u)v$ for $1\leq i\leq m+n$.
\end{enumerate}

\bth[\cite{Zhang 1996}]
$L(\la(u))$ is finite dimensional if and only if its highest weight $\la(u)$ satisfies the following conditions:
\[
\frac{\la_k(u)}{\la_{k+1}(u)}=\frac{P_k(u+(-1)^{\bk})}{P_k(u)},
\]
\[
\frac{\la_m(u)}{\la_{m+1}(u)}=\frac{Q_0(u)}{Q_1(u)},
\]
for $1\leq k<m+n$ and $k\neq m$, where $P_k(u)$ is a polynomial in $u$ and
\[
Q_0(u)=\prod_{s=1}^{N}(u+a_i),\quad Q_1(u)=\prod_{s=1}^{N}(u-b_i),
\]
such that $Q_0(u)$ and $Q_1(u)$ are coprime.
\eth
In fact, up to an automorphism $\om_f$ \eqref{omf} of $\Y(\gl_{m|n})$, Drinfeld polynomials  $P_k$, $Q_i$  uniquely determine a finite dimensional simple $\Y(\gl_{m|n})$-module.

\subsection{Gelfand-Tsetlin character}
Given $\La\in\mathscr{S}_{\la}$, define
\[
\chil(u)=\prod_{k=1}^{i}\frac{u+l_{i,k}}{u-k+1},
\]
and
\[
\chjl(u)=\prod_{k=1}^{m}\frac{u+l_{m+j,k}}{u-k+1}\prod_{k'=1}^{j}\frac{u-m+k'}{u+l_{m+j,m+k'}},
\]
for $1\leq i\leq m$ and $1\leq j\leq n$.
\ble\label{lem:chkL}
Let $\La$ be a Gelfand-Tsetlin pattern in $\mathscr{S}_{\la}$ and let $\xiL$ be the corresponding Gelfand-Tsetlin basis vector.
Then we have
\begin{enumerate}
\item $B_i(u)\xiL=\chil(u)\xiL$,
\item $B_{m+j}(u)\xiL=\chjl(u)\xiL$.
\end{enumerate}
\ele

\subsection{Nonsuper case}
We recall the following irreducibility criterion of $L(\la)\ot L(\ga)$ with both $\la$ and $\mu$  integral dominant, given by Molev in \cite{Mol 2002} for the Yangian $\Y(\gl_n)$.

Denote $\la_i-i+1$ and $\ga_i-i+1$ by $l_i$ and $r_i$ respectively for all $1\leq i\leq n$, and
 define the set $\langle l_j,l_i\rangle$ by
\beq
\langle l_j,l_i\rangle=\{l_j,l_j+1,l_j+2,\ldots,l_i\}\backslash\{l_j,l_{j-1},\ldots,l_i\},
\eeq
for $1\leq i< j\leq n$. The sets $\langle r_j,r_i\rangle$ are defined in the same way.

Then the irreducibility criterion is described by
\bth[\cite{Mol 2002}]\label{thm:nonsuper}
Let $\la$ and $\ga$ be both integral dominant weight of $\gl_{n}$.
Then $L(\la)\ot L(\ga)$ is a simple $\Y(\gl_{n})$-module if and only if for each pair of indices $1\leq i< j\leq n$ holds
\beql{ns noncro}
r_j,r_i\notin \langle l_j,l_i\rangle\quad\text{ or }\quad l_j,l_i\notin \langle r_j,r_i\rangle.
\eeq
\eth

\section{Irreducibility criterion}\label{sec:cri}
Recall that for a $\gl_{m|n}$-weight $\la=(\la_1,\ldots,\la_{m+n})$, the shift by $h\in \BC$ is defined as
\[
\la+h=(\la_1+h,\ldots,\la_m+h,\la_{m+1}-h,\ldots,\la_{m+n}-h).
\]

\bpr\label{prop:shif}
The $\Y(\gl_{m|n})$-module
\[
L_g(\la)\ot L_h(\ga)
\]
is simple if and only if the $\Y(\gl_{m|n})$-module
\[
L(\la+g)\ot L(\ga+h)
\]
is simple.
\epr
\bpf
This proof is analogous to the one in \cite{Mol 2002}.
Suppose that $L_g(\la)\ot L_h(\ga)$ is simple.
Let $\xi_0$ and $\eta_0$ denote  highest weight vectors of $\gl_{m|n}$-modules $L(\la)$ and $L(\ga)$ respectively.
It is not difficult to see that $\xi_0\ot\eta_0$ is a highest weight vector of $L_g(\la)\ot L_h(\ga)$ with highest weight $(\la_1(u),\ldots,\la_{m+n}(u))$, where
\[
\la_i(u)=(1+(-1)^{\bi}\frac{\la_i}{u+g})(1+(-1)^{\bi}\frac{\ga_i}{u+h}).
\]
Consider the automorphism $\om_f:t_{ij}(u)\mapsto f(u)t_{ij}(u)$ of $\Y(\gl_{m|n})$ for$f(u)=(1+gu^{-1})(1+hu^{-1})$ (cf. Section 3, \eqref{omf}).
Applying the automorphism $\om_f$ to $L_g(\la)\ot L_h(\ga)$ we obtain a simple $\Y(\gl_{m|n})$-module $\tilde{L}$ with highest weight $(\tilde{\la}_1(u),\ldots,\tilde{\la}_{m+n}(u))$ where
\[
\tilde{\la}_i(u)=(1+\frac{(-1)^{\bi}(\la_i+(-1)^{\bi}g)}{u})(1+\frac{(-1)^{\bi}(\ga_i+(-1)^{\bi}h)}{u}).
\]
Consider the subquotient of $L(\la+g)\ot L(\ga+h)$ generated by a highest weight vector $\xi_0'\ot\eta_0'$.
Since
\[
t_{ii}(u)\xi_0'\ot \eta_0'=\tilde{\la}_i(u)\xi_0'\ot \eta_0',
\]
then $\tilde{L}$ is isomorphic to the subquotient of $L(\la+g)\ot L(\ga+h)$. However, $\tilde{L}$ and $L(\la+g)\ot L(\ga+h)$ have the same dimension, and hence they are isomorphic. It implies that $L(\la+g)\ot L(\ga+h)$ is simple. The converse statement is verified by reversing the argument.
\epf

The above proposition allows us to reduce the problem of determining the irreducibility of $L_g(\la)\ot L_h(\ga)$, where each of $\la$ and $\ga$ is either essentially typical or covariant to that of $L(\la+g)\ot L(\ga+h)$, where $\la+g$ and $\ga+h$ are either essentially typical or shifted covariant.

Therefore, throughout the remainder of the paper, we shall restrict our attention to tensor products $L(\la)\ot L(\ga)$, where each of $\la$ and $\ga$ is either essentially typical or shifted covariant.

\subsection{Sufficient condition}

Let $\xi_0$ and $\eta_0$ denote the highest vectors of the $\gl_{m|n}$-modules
$L(\la)$ and $L(\ga)$, respectively.
Let $N$ be a nonzero $\Y(\gl_{m|n})$-submodule of $L(\la)\ot L(\ga)$.
It is clear that $N$ must
contain a singular vector $\zeta$.
The key part of the proof of the theorem is to show that
\beql{sing}
\zeta={\rm const}\cdot\xi_0\ot\eta_0.
\eeq
Then considering dual modules we show that
the vector $\xi_0\ot\eta_0$ is a cyclic generator of $L(\la)\ot L(\ga)$.

In the remainder of this paper, we replace $t_{ij}(u)$ by $u^2t_{ij}(u)$ in order to avoid the appearance of denominators.
For convenience, we continue to denote it by $t_{ij}(u)$.
Using the Gelfand-Tsetlin basis, the singular vector
$\zeta$ is uniquely written in the form
\beql{zeta}
\zeta=\sum_{\substack{\La\in\mathscr{S}_{\la}\\ \Ga\in\mathscr{S}_{\ga}}}c_{\La,\Ga}\xiL\ot \eta_{\Ga},
\eeq
where $\xi_{\La}$ and $\eta_{\Ga}$ are Gelfand-Tsetlin basis elements of $L(\la)$ and $L(\ga)$ parameterized by  admissable patterns $\La$ and $\Ga$ respectively, and $c_{\La,\Ga}\in \BC$. We say a pattern $\La\in\mathscr{S}_{\la}$ occurs in $\zeta$ if there exists some $\Ga\in\mathscr{S}_{\ga}$ such that $c_{\La,\Ga}\neq 0$.

For the diagonal Cartan subalgebra $\mathfrak h$ of $\gl_{m|n}$,
we denote by $\ve_i$  the basis vector	of $\mathfrak h^*$
dual to the element $E_{ii}$ for $1\leq i\leq m+n$, so that the $m+n$-tuple $\la$ can be identified
with the element $\la_1\ve_1+\cdots+\la_{m+n}\ve_{m+n}\in \mathfrak h^*$.
We shall use a standard partial ordering on the weights of $L(\la)$.
Given two weights $v,w\in\mathfrak h^*$,
we write $v\preceq w$ if $w-v$ is a
$\ZZ_+$-linear combination
of the simple roots
$\ve_a-\ve_{a+1}$. Equivalently, $v\preceq w$ if and only if
\beql{eqp}
w-v=\sum_{a=1}^{m+n}p_a\ve_a,
\eeq
with the conditions
\beql{eqpco}
p_1,\ p_1+p_2,\ \dots,\ p_1+\cdots+p_{m+n-1}\in\ZZ_+,\quad
p_1+\cdots+p_{m+n}=0.
\non
\eeq

The embedding
\beql{emb}
\U(\gl_{m|n})\hookrightarrow\Y(\gl_{m|n}),\qquad E_{ij}\mapsto t_{ij}^{(1)}(-1)^{\pa{i}},
\eeq
defines the natural
$\U(\gl_{m|n})$-module  structure
on $L(\la)\ot L(\ga)$. We shall usually identify the operators
$E_{ij}$ and $t_{ij}^{(1)}(-1)^{\pa{i}}$.
The vector $\zeta$ is clearly a $\gl_{m|n}$-singular vector. In particular,
it is a weight vector. Since the basis $\{\xiL\}$ consists of
weight vectors, each element $\etaL$ occurring in \eqref{zeta} is also
a $\gl_{m|n}$-weight vector. Moreover, all elements $\xiL\ot\eta_{\Ga}$ in \eqref{zeta} have
the same $\gl_{m|n}$-weight.

We shall denote the weight of the vector $\xiL$, or the weight
of the pattern $\La$, by $w(\La)$.
It can be deduced from \eqref{action h} that
\beql{weigt}
w(\La)=w_1\ve_1+\cdots+w_{m+n}\ve_{m+n},\qquad w_k=
\sum_{i=1}^k\la_{ki}-\sum_{i=1}^{k-1}\la_{k-1,i}.
\eeq

Consider the set of patterns occurring in $\zeta$
and suppose that $\La^0=(\la_{ij}^0)$ is a minimal element of this set with respect to the
partial ordering on the weights $w(\La)$. In other words, if $\La$ occurs in $\zeta$
and $w(\La)\preceq w(\La^0)$
then $w(\La)=w(\La^0)$.

\ble\label{lem:ga0}
$a_{\La^0,\Ga}$ equals zero unless $\eta_{\Ga}$ is the highest weight vector in $L(\ga)$.
\ele
\bpf
We have $t_{ii+1}(u)\ts\zeta=0$ for $i=1,\dots,m+n-1$. Therefore, by
\eqref{copro1},
\beql{bzet}
\sum_{j=1}^{m+n}\sum_{\substack{\La\in\mathscr{S}_{\la}\\ \Ga\in\mathscr{S}_{\ga}}}c_{\La,\Ga}t_{ij}(u)\xiL\ot t_{ji+1}(u)\eta_{\Ga}=0.
\eeq
We compare the coefficients of $\xiLo$ on the both sides of the equation. It is clear that the indeterminant $u$  only occurs in the summands when $j=i$ or $j=i+1$. Since $w(\La^0)$ is minimal, $t_{ii+1}(u)\ts\xiL$ cannot be a multiple of $\xiLo$. Thus the coefficient of $\xiLo$ containing $u$ in a summand can only occur in the case  $j=i$. It implies that $t_{ii+1}(u)\eta_{\Ga}=0$ if $c_{\La^0,\Ga}\neq 0$. Therefore, the corresponding $\eta_{\Ga}$ is a highest weight vector in $L(\ga)$.
\epf

\ble\label{lem:uni}
The pattern $\La^0$ is determined uniquely.
Moreover, if a pattern $\La$ occurs in $\zeta$
then $w(\La)\succeq w(\La^0)$.
\ele
\bpf

Since $\zeta$ is a highest weight vector, we have
\[
B_k(u)\zeta=f_k(u)\zeta,
\]
where $f_k(u)\in \BC(u)$ for $1\leq k\leq m+n$.

The coproduct of $B_k(u)$ can be written as
\[
\Delta(B_k(u))=\sum B_k(u)_{(1)}\otimes B_k(u)_{(2)}.
\]
Assume that there exists another pattern $\La'$ such that $w(\La')$ is  minimal.
Consider
\[
B_k(u)\zeta=B_k(u)\sum_{\substack{\La\in\mathscr{S}_{\la}\\ \Ga\in\mathscr{S}_{\ga}}}c_{\La,\Ga}\xiL\ot \eta_{\Ga}=\sum_{\substack{\La\in\mathscr{S}_{\la}\\ \Ga\in\mathscr{S}_{\ga}}}c_{\La,\Ga}B_k(u)_{(1)}\xiL\ot B_k(u)_{(2)}\eta_{\Ga}.
\]

Due to \eqref{com1}, for $k\leq m$, $B_k(u)_{(1)}$ equals the coefficient of $E_{1a_1}\ot E_{2a_2}\ot\cdots\ot E_{ka_k}$ in
\[
A_kT_{1[1]}(u)\cdots T_{k[1]}(u-k+1)
\]
for some $1\leq a_1,\ldots,a_k\leq m+n$. We denote such coefficient by $t_{a_1\ldots a_k}^{1 \ldots k}(u)$.
It is clear that $w(t_{a_1\ldots a_k}^{1 \ldots k}(u)\xi)\geq w(\xi)$ for any weight vector $\xi$.
Hence $f_k(u)$ must be equal to the coefficient of $B_k(u)\xi_{\La^0}$. But it also equals the coefficient of $B_k(u)\xi_{\La'}$. Now it follows from Lemma \ref{lem:chkL} that $\la^0_{ij}=\la'_{ij}$ for $1\leq j\leq i\leq m$.

Similarly, thanks to \eqref{com2}, for $k>m$, $B_k(u)_{(1)}$ is the coefficient of $E_{1a_1}\otimes\cdots\otimes E_{ka_k}$ in
\[
A_mT_{1[1]}(u)\cdots T_{m[1]}(u-m+1)H_{k-m}^{\circ}((T(u)^*)^{\mathbb{K}}_{\mathbb{M}+\mathbb{N}})_{m+1[1]}\cdots ((T(u)^*)^{\mathbb{K}}_{\mathbb{M}+\mathbb{N}})_{k[1]}
\]
for some $1\leq a_1,\ldots,a_k\leq m+n$. Here, $\mathbb{K}$ is the set $\{1,2,\ldots,k\}$. We also denote such coefficient by $t_{a_1\ldots a_k}^{1 \ldots k}(u)$.

If  $\la^0_{kj}\neq \la'_{kj}$ for some $k>m$ and $1\leq j\leq k$ then the coefficients of $B_k(u)\xi_{\La^0}$ and $B_k(u)\xi_{\La'}$ cannot be equal by Lemma \ref{lem:chkL}. Moreover, the coefficients are rational function of the form $\frac{G(u)}{H(u)}$ and $\frac{G'(u)}{H'(u)}$ respectively, with $G(u)=G'(u)$ and $\text{deg}(H(u))=\text{deg}(H'(u))$.
If there exists some $\La_a$ occurring in $\zeta$ such that $t_{a_1\ldots a_k}^{1 \ldots k}(u)\xi_{\La_a}=\frac{G_a(u)}{H_a(u)}\xi_{\La^0}$, then $\text{deg}(H_a(u))<\text{deg}(H(u))$. By expanding them in $\BC((u^{-1}))$, we have that
$\xi_{\La}\ot\eta_{\Ga}$ and $\xi_{\La'}\ot\eta_{\Ga}$ have  distinct coefficients in $B_k(u)\zeta$.
Consequently, $\la^0_{ij}=\la'_{ij}$ for all $1\leq j\leq i\leq m+n$.

\epf

\ble\label{lem:luc}
If $0\leq i_0\leq n-1$ is the maximal such that $\la_{m+i_0,k}^{0}<\la_{m+i_0+1,k}^{0}$ for any $1\leq k\leq m$, then
\[
\la_{m,k}^{0}<\la_{m+1,k}^{0}<\cdots<\la_{m+i_0,k}^{0}<\la_{m+i_0+1,k}^{0}.
\]
\ele

\bpf
We write \eqref{zeta} as
\beql{zetai}
\zeta=\sum_{p=0}^Kc_p\xi_{\La^p}\ot\eta_{\Ga^p},
\eeq
where $\La^0$ is the minimal element among $\{\La^0,\La^1,\ldots,\La^K\}$ with respect to the
partial ordering on the weights $w(\La)$.
According the condition on $i_0$ and Theorem \ref{thm:GT basis}, we have that
\[
\La'=\La^0+\delta_{m,k}+\delta_{m+1,k}+\cdots+\delta_{m+i_0,k}
\]
is an admissable pattern.
Due to Theorem \ref{thm:GT basis}, $t_{m,m+i_0+1}(u)\xi_{\La^0}=C_0\xi_{\La'}$ for some nonzero $C_0\in\BC$. But
\beql{elm}
0=t_{m,m+i_0+1}(u)\zeta=\sum_p\sum_qc_pt_{m,q}(u)\xi_{\La^p}\ot t_{q,m+i_0+1}(u)\eta_{\Ga^p}.
\eeq
Since $t_{m,m+i_0+1}(u)\xi_{\La^0}\ot t_{m+i_0+1,m+i_0+1}(u)\eta_{\Ga^0}=f(u)\xi_{\La'}\ot \eta_{\Ga^0}$ for some polynomial $f(u)$ in $u$ with degree $1$, there must exist a summand $c't_{m,m}(u)\xi_{\La'}\otimes t_{m,m+i_0+1}(u)\eta_{\Ga'}$ in \eqref{elm} to eliminate $\xi_{\La'}\ot \eta_{\Ga^0}$.
Thus, $\La'$ occurs in $\zeta$.
Without loss of generality, we assume that $\La^1=\La'$.
The weight of $\Ga^1$ equals $w(\Ga^0)-(\ve_m-\ve_{m+i})$. Thus
\[
\Ga^1=\Ga^0-\delta_{m,a_0}-\delta_{m+1,a_1}-\cdots-\delta_{m+i_0,a_{i_0}}
\]
where $1\leq a_j\leq m+j$ for $0\leq j\leq i_0$.

It is easy to see that $\Ga^2=\Ga^1+\delta_{m+i_0,a_{i_0}}$ is an admissable pattern.
Thus, $t_{m+i_0,m+i_0+1}(u)\eta_{\Ga^1}=C_1\eta_{\Ga^2}$ for some nonzero $C_1$.
By the same argument,  we get that $\La''=\La^1-\delta_{m+i_0,a}$ occurs in $\zeta$ for some $1\leq a\leq m+i_0$.
Lemma \ref{lem:uni} implies that $a$ must be equal to $k$.
In fact, if it is not this case then we can apply $t_{i,i+1}(u)$ to $\eta_{\Ga^2}$ repeatedly to get $\eta_{\Ga^0}$, and similarly to the argument for the occurrence of $\La^1$, we  show that a pattern $\La^o$  occurs in $\zeta$.
Since $a\neq k$, we have that  $\La^o\neq \La^0$ but $w(\La^0)=w(\La^o)$, which contradicts to Lemma \ref{lem:uni}.

Consequently, $\La''=\La^0+\delta_{m,k}+\delta_{m+1,k}+\cdots+\delta_{m+i_0-1,k}$ appears in $\zeta$,
whence $\la^0_{m+i_0,k}>\la^0_{m+i_0-1,k}$.
The proof is then completed by repeating the above argument.

\epf

\ble\label{lem:wadd con}
 If
$-\ga_{m+1}\neq \la_m$ then $\la_{m+i,k}^0=\la_k$ for all $0\leq i\leq n-1$ and $1\leq k\leq m$.
\ele
\bpf
Suppose that there exists some $0\leq i\leq n-1$ such that $\la^0_{m+i,m}<\la_m$. Let $i_0$ be the maximal such that $\la_{m+i_0,m}^0<\la_{m+i_0+1,m}^0$ and $\La^1$ be the pattern described in the proof of Lemma \ref{lem:luc}.
According to Theorem \ref{thm:GT basis} and Lemma \ref{lem:luc}, we have $\la^0_{m+j,m}=\la^0_{m+j+1,m}-1$ for $0\leq j\leq i_0$.
Thus, in $\La^1$,
\[
\la_{m,m}^{1}<\la_{m+1,m}^{1}<\cdots<\la_{m+i_0,m}^{1}=\la_{m+i_0+1,m}^{1}.
\]

Applying $t_{m,m+i_0}(u)$ to $\zeta$ and using the same argument in the proof of Lemma \ref{lem:luc} we have that $\La^2=\La^1+\delta_{m,m}+\delta_{m+1,m}+\cdots+\delta_{m+i_0-1,m}$ occurs in $\zeta$. Repeating this process, we get that $\La^{i_0}$ occurs in $\zeta$ and
\[
\la_{m,m}^{i_0}<\la_{m+1,m}^{i_0}=\la_{m+2,m}^{i_0}.
\]

Applying $t_{m,m+1}(u)$ to $\zeta$, we get that
\[
0=t_{m,m+1}(u)\zeta=\sum_{p=1}^K\sum_{q=1}^{m+n}c_pt_{mq}(u)\xi_{\La^p}\ot t_{q,m+1}(u)\eta_{\Ga^p}.
\]
Tensor products of Gelfand-Tsetlin patterns form a basis of the tensor product of modules. In this basis we have

\[
\begin{split}
t_{m,m+1}(u)\xi_{\La^{i_0}}\ot t_{m+1,m+1}(u)\xi_{\Ga^{i_0}}=&-E_{m,m+1}\xi_{\La^{i_0}}\ot (u-E_{m+1,m+1})\eta_{\Ga^{i_0}}\\
=&C_2\xi_{\La^{i_0}+\delta_{mm}}\ot (u-\ga_{m+1})\eta_{\Ga^{i_0}}+\text{other terms}.
\end{split}
\]
It follows from Theorem \ref{thm:GT basis}  that $C_2\neq 0$.
To eliminate this term, we consider $t_{mm}(u)\xi_{\La^{i_0}+\delta_{mm}}\ot t_{m,m+1}(u)\eta_{\Ga^{i_0}-\delta_{mm}}$.
Without loss of generality we can assume that $\la^{i_0}_{ij}=\la^{i_0}_{mj}$ for all $1\leq j\leq i\leq m$ (otherwise replace $\zeta$ using the action of $\Y(\gl_{mn})$  and the argument as in the proof of Lemma \ref{lem:luc}). But
\[
\begin{split}
t_{mm}(u)\xi_{\La^{i_0}+\delta_{mm}}\ot t_{m,m+1}(u)\eta_{\Ga^{i_0}-\delta_{mm}}=&(u+E_{mm})\xi_{\La^{i_0}+\delta_{mm}}\ot (-E_{m,m+1})\eta_{\Ga^{i_0}-\delta_{mm}}\\
=&(u+\la^{i_0}_{mm}+1)\xi_{\La^{i_0}+\delta_{mm}}\ot C_3\eta_{\Ga^{i_0}}+\text{other terms}
\end{split}
\]
for some nonzero $C_3$.

Consequently,  we have $C_2'(u-\ga_{m+1})+C_3'(u+\la^{i_0}_{mm}+1)=0$ for some nonzero constants $C_2'$ and $C_3'$. Then $-\ga_{m+1}=\la^{i_0}_{mm}+1=\la_m$, which contradicts to the condition $-\ga_{m+1}\neq \la_m$. This completes the statement for $k=m$. The case $k<m$ can be shown similarly.

\epf
Recall that
\beq
l_{i}=\la_{i}-i+1, (1\leq i\leq m); \quad l_{j}=-\la_{j}+j-2m, (m+1\leq j\leq m+n)
\eeq
and set
\beq
r_{i}=\ga_{i}-i+1, (1\leq i\leq m); \quad r_{j}=-\ga_{j}+j-2m, (m+1\leq j\leq m+n).
\eeq
 For pairs of indices $(i,j)$ and $(p,q)$ such that
$1\leq i<j\leq m$ and  $m+1\leq p<q\leq m+n$, we define
\beq
\langle l_j,l_i\rangle=\{l_j,l_j+1,l_j+2,\ldots, l_i\}\backslash\{l_j,l_{j-1},\ldots, l_i\}.
\eeq
and
\beq
\langle l_p,l_q\rangle=\{l_p,l_p+1,l_p+2,\ldots,l_q\}\backslash\{l_p,l_{p+1},\ldots, l_q\}.
\eeq
Similarly, we define $\langle r_j,r_i\rangle$ and $\langle r_p,r_q\rangle$.

\ble\label{lem:la0}
Suppose that $\la$ and $\ga$ satisfy the following conditions:
\[
l_j,l_i\notin \langle r_j,r_i\rangle,\quad l_p,l_q\notin\langle r_p,r_q\rangle \text{ and }-\ga_{m+1}\neq \la_m,
\]
for $1\leq i<j\leq m$ and $m+1\leq p<q\leq m+n$.
Then $\la_{ij}^0=\la_j$ for all $1\leq i\leq m+n-1$ and $1\leq j\leq i$.
\ele
\bpf
Since $-\ga_{m+1}\neq \la_m$, then $\la_{mk}^0=\la_k$ for all $0\leq i\leq n-1$ and $1\leq k\leq m$.
Regarding $\zeta$ as a $\Y(\gl_m)$ highest weight vector,  we see that $\la^0_{ij}=\la_{mj}^0=\la_j$ for $1\leq j\leq i\leq m$ by  Theorem \ref{thm:nonsuper}.
Similarly,  $\la^0_{ij}=\la_j$ for $m<j\leq i<m+n$ by regarding $\zeta$ as a $\Y(\gl_{0|n})$ highest weight vector.
\epf

We have proved that $\zeta=\text{const}\cdot\xi_0\ot \eta_0$. Using the same argument, we can see that if
\[
r_j,r_i\notin \langle l_j,l_i\rangle,\ r_p,r_q\notin\langle l_p,l_q\rangle \text{ and } -\ga_{m+1}\neq \la_m
\]
for $1\leq i<j\leq m$ and $m+1\leq p<q\leq m+n$, then again $\zeta=\text{const}\cdot\xi_0\ot \eta_0$.

Furthermore, if $\la$ and $\ga$ satisfy the  conditions
\begin{enumerate}
\item $l_j,l_i\notin \langle r_j,r_i\rangle$ or $r_j,r_i\notin \langle l_j,l_i\rangle$,
\item $l_p,l_q\notin\langle r_p,r_q\rangle$ or $r_p,r_q\notin\langle l_p,l_q\rangle$,
\item $-\ga_{m+1}\neq \la_m$,
\end{enumerate}
then any singular vector in $L(\ga)\ot L(\la)$ is of the form $\text{const}\cdot \eta_0\ot \xi_0$.

\bde\label{def:noncro}
We say $\la$ and $\ga$ satisfy the \emph{super non-crossing condition} if for all $1\leq i<j\leq m$ and $m+1\leq p<q\leq m+n$, the following conditions hold:
\begin{enumerate}
\item $l_j,l_i\notin \langle r_j,r_i\rangle$ or $r_j,r_i\notin \langle l_j,l_i\rangle$;
\item $l_p,l_q\notin\langle r_p,r_q\rangle$ or $r_p,r_q\notin\langle l_p,l_q\rangle$;
\item  $-\la_{m+1}\neq\ga_m$ and $-\ga_{m+1}\neq \la_m$.
\end{enumerate}
\ede

There is an equivalent definition given by Molev \cite{Mol 2002}, which perhaps better reflects the meaning of ``non-crossing''.
Two disjoint finite subsets $A$ and $B$ of $\mathbb{Z}$ are called crossing if there exist element $a_1,a_2\in A$ and $b_1,b_2\in B$ such that
\[
a_1<b_1<a_2<b_2\quad\text{ or }\quad b_1<a_1<b_2<a_2.
\]
Otherwise, $A$ and $B$ are called non-crossing.

For any $\gl_{m|n}$-weight $\la$
introduce the sets
\[
\mathscr{E}_{\la}=\{\la_1,\la_2-1,\ldots,\la_m-m+1\}
\]
and
\[
\mathscr{O}_{\la}=\{-\la_{m+1}-m+1,-\la_{m+2}-m+2,\ldots,-\la_{m+n}-m+n\}.
\]

Then the super non-crossing condition for weights $\lambda$ and $\gamma$
 can be written equivalently as follows:

\begingroup
\begin{enumerate}
\item $\mathscr{E}_{\la}\backslash \mathscr{E}_{\ga}$ and $\mathscr{E}_{\ga}\backslash \mathscr{E}_{\la}$ are non-crossing;
\item $\mathscr{O}_{\la}\backslash \mathscr{O}_{\ga}$ and $\mathscr{O}_{\ga}\backslash \mathscr{O}_{\la}$ are non-crossing;
\item $-\la_{m+1}\neq \ga_m$ and $-\ga_{m+1}\neq \la_m$.
\end{enumerate}
\addtocounter{thm}{-1}
\endgroup

\bth\label{thm:scon}
Let $\la$ and $\ga$ be shifted covariant or essentially typical $\gl_{m|n}$-weights.
If they satisfy the super non-crossing condition, then $V=L(\la)\ot L(\ga)$ is a simple $\Y(\gl_{m|n})$-module.
\eth
\bpf
We first prove that if
\begin{enumerate}
\item  $-\la_{m+1}\neq\ga_m$ and $-\ga_{m+1}\neq \la_m$;
\item  $l_j,l_i\notin \langle r_j,r_i\rangle$, $l_p,l_q\notin\langle r_p,r_q\rangle$,
\end{enumerate}
then $V$ is simple.

By Lemma \ref{lem:la0}, it suffices to show that $V$ is generated by $\zeta$.
Let $N=\Y(\gl_{m|n})\zeta$. It follows from Proposition \ref{pro:dual} that  the dual module $V^*$
is isomorphic to $\overline{L(\la)\ot L(\ga)}$, and hence to $L(\ga)\ot L(\la)$ through the map $\xi\ot \eta\mapsto \eta\ot\xi$.

If $N$ is a proper submodule of $V$ then $\text{Ann }{N}$ is a proper submodule of $V^*$ which contains a nonzero singular vector of $V^*$. Since $V^*$ is isomorphic to $L(\ga)\ot L(\la)$, this singular vector must be of the form $f_{\eta_0}\ot f_{\xi_0}\in \text{Ann }{N}$, where $f_{\xi}$ is the dual basis with respect to $\xi$ (according to the same argument as in the proof  that  singular vectors in $L(\la)\ot L(\ga)$ are all of the form $\rm{const}\cdot \xi_0\ot \eta_0$).  This gives a contradiction since $(f_{\xi_0}\ot f_{\eta_0})(\xi_0\ot \eta_0)\neq 0$.   Hence $N=V$ proving the irreducibility of $V$.

For the same reasons, all other circumstances of the super non-crossing condition also imply the irreducibility of $L(\la)\ot L(\ga)$.
\epf

\subsection{Necessary condition}
In this section our aim is to show that the sufficient condition given in Theorem \ref{thm:scon} is also necessary.

For $\la=(\la_1,\la_2,\ldots,\la_m,\la_{m+1},\ldots,\la_{m+n})$ let $\la^e=(\la_1,\ldots,\la_m)$ and $\la^o=(\la_{m+1},\ldots,\la_{m+n})$ be the even and the odd part of $\la$ respectively.

\ble
Let $\la$ and $\ga$ be shifted covariant or essentially typical $\gl_{m|n}$-weights. Suppose that the $Y(\gl_{m|n})$-module $L(\la)\ot L(\ga)$  is simple. Then the $Y(\gl_{m})$-module  $L(\la^e)\ot L(\ga^e)$ and the $Y(\gl_{n})$-module
  $L(\la^o)\ot L(\ga^o)$ are simple.
\ele
\bpf
We will identify $L(\la^e)$ and $L(\ga^e)$ with the $\U(\gl_m)$-spans of the highest weight vectors $\xi_0$  and $\eta_0$ in  $L(\la)$ and in $L(\ga)$ respectively.
The generators $E_{i,m+j}$ of $\gl_{m|n}$ with $1\leq i< m+j\leq m+n$ annihilate $L(\la^e)$ and $L(\ga^e)$.
Hence, the subspace $L(\la^e)\ot L(\ga^e)$ of $L(\la)\ot L(\ga)$ is invariant with respect to the action of the subalgebra $\Y(\gl_m)$ of $\Y(\gl_{m|n})$.
Moreover, this  $\Y(\gl_m)$-module structure can be obtained by regarding $L(\la^e)$ and $L(\ga^e)$ as evaluation $\Y(\gl_m)$-modules  and applying the coproduct in $\Y(\gl_m)$.

Suppose that there is a nonzero submodule of $L(\la^e)\ot L(\ga^e)$ which does not contain the vector $\xi_0\ot \eta_0$.
Then this submodule contains a nonzero vector $\zeta'$ annihilated by all $t_{i,i+1}(u)$ with $1\leq i\leq m-1$. However, we also have $t_{m+j,m+j+1}(u)\zeta'=0$ for any $0\leq j\leq n-1$ which follows from \eqref{copro1}.
This implies that $\zeta'$ is  a singular vector in $L(\la)\ot L(\ga)$ which shows that $L(\la)\ot L(\ga)$ is not simple. It contradicts to our assumption.

Suppose now that the $\Y(\gl_m)$-submodule of $L(\la^e)\ot L(\ga^e)$ generated by $\xi_0\ot \eta_0$ is proper.
Equip the dual space $(L(\la^e)\ot L(\ga^e))^*$ with a $\Y(\gl_m)$-module structure through the antiautomorphism $\st$ given in Proposition \ref{pro:anti-auto} in the case $n=0$.
We have shown  in Proposition \ref{pro:dual} that $(L(\la^e)\ot L(\ga^e))^*$ is isomorphic to $L(\ga^e)\ot L(\la^e)$.
Since $\xi_0\ot \eta_0$ generate a proper submodule of $L(\ga^e)\ot L(\la^e)$, we know (cf. the argument before Proposition \ref{pro:dual}) that its annihilator in $(L(\la^e)\ot L(\ga^e))^*\cong L(\ga^e)\ot L(\la^e)$ is  a nonzero submodule which does not contain $\eta_0\ot \xi_0$.

On the other hand, since $L(\la)\ot L(\ga)$ and the simple subquotient of $L(\ga)\ot L(\la)$ generated by the tensor product of highest weight vectors share the same Drinfeld polynomials, $L(\ga)\ot L(\la)$ is simple \cite{Zhang 1996}. Similarly to the case of $L(\la^e)\ot L(\ga^e)$, the proper submodule of $L(\ga^e)\ot L(\la^e)$ must contain $\eta_0\ot \xi_0$, which leads to a contradiction.

The second statement follows by a similar argument identifying $L(\la^o)$ and $L(\ga^o)$ with the $\U(\gl_{0|n})$-spans of the highest weight vectors $\xi_0$ and $\eta_0$ in $L(\la)$  and $L(\ga)$ respectively.
\epf

According to \cite{Mol 2002}, the irreducibility of $L(\la)\ot L(\ga)$ implies that  the first two conditions in Definition \ref{def:noncro} hold.
If the third condition fails, for example, $\la_m=-\ga_{m+1}$, we construct a nonzero proper submodule of $L(\la)\ot L(\ga)$.

In fact, this instance can be found in the proof of Lemma \ref{lem:wadd con}.
Let  $$\zeta'=c_1\xi_{\La^1}\ot \eta_{\Ga^1}+c_2\xi_{\La^2}\ot \eta_{\Ga^2}$$ for some complex numbers $c_1,c_2$, where $\La^1$ and $\Ga^2$ are admissable patterns with
\[
\la^1_{ij}=\begin{cases}
\la_m-1,&\text{ if } (i,j)=(m,m),\\
\la_j,&\text{otherwise},
\end{cases}
\]
and
\[
\ga^2_{ij}=\begin{cases}
\ga_m-1,&\text{ if } (i,j)=(m,m),\\
\ga_j,&\text{otherwise},
\end{cases}
\]
while $\Ga^1$ and $\La^2$ are admissable patterns corresponding to highest weight vectors in $L(\ga)$ and $L(\la)$ respectively.

Now it is clear that $t_{i,i+1}(u)\zeta'=0$ unless $i=m$.
Moreover,
\[
\begin{split}
t_{m,m+1}(u)\zeta'=&c_1t_{m,m+1}(u)\xi_{\La^1}\ot t_{m+1,m+1}(u)\eta_{\Ga^1}+c_2t_{mm}(u)\xi_{\La^2}\ot t_{m,m+1}(u)\eta_{\Ga^2}\\
=&c_1E_{m,m+1}\xi_{\La^1}\ot (u-E_{m+1,m+1})\eta_{\Ga^1}+c_2(u+E_{mm})\xi_{\La^2}\ot E_{m,m+1}\eta_{\Ga^2}\\
=&c_1C_1(u-\ga_{m+1})\xi_{\La^2}\ot \eta_{\Ga^1}+c_2C_2(u+\la_m)\xi_{\La^2}\ot \eta_{\Ga^1},
\end{split}
\]
where $C_1,C_2$ are nonzero complex numbers.
Choosing $c_1=\frac{1}{C_1}$ and $c_2=-\frac{1}{C_2}$,
 we get that $t_{m,m+1}(u)\zeta'=0$. Thus, $\zeta'$ is a singular vector not proportional to $\zeta$ and the submodule of $\Y(\gl_{m|n})$ generated by $\zeta'$ is proper.

As a result, we have proved the following theorem.
\bth\label{thm:ncon}
Let $\la$ and $\ga$ be shifted covariant or essentially typical $\gl_{m|n}$-weights. Suppose that the $Y(\gl_{m|n})$-module $L(\la)\ot L(\ga)$ is simple. Then  $\la$ and $\ga$ satisfy the super non-crossing condition.
\eth

Theorem \ref{thm:ncon} completes the proof of our main result (Theorem A) from the introduction which gives the irreducibility criterion for the tensor product $L(\la)\ot L(\ga)$.

\bre\label{coro:multi}
In \cite{NT 2002}, the binary property for tensor products of skew modules over the Yangian $\Y(\gl_n)$ was established. More precisely, the tensor product
\[
L_{h_1}(\la^{(1)}\slash\mu^{(1)})\otimes
L_{h_2}(\la^{(2)}\slash\mu^{(2)})\otimes
\cdots\otimes
L_{h_k}(\la^{(k)}\slash\mu^{(k)})
\]
is simple if and only if, for every pair $(i, j)$ with $1\le i<j\le k$, the tensor product
\[
L_{h_i}(\la^{(i)}\slash\mu^{(i)})
\otimes
L_{h_j}(\la^{(j)}\slash\mu^{(j)})
\]
is simple.
Here, $L_h(\la\slash\mu)$ denotes the skew module of $\Y(\gl_n)$ associated with the skew Young diagram $\la\slash\mu$. When $\mu$ degenerates to $0$, the skew module reduces to the evaluation module $L_h(\la)$.

One can observe that the same method applies to the covariant weight modules of super Yangian $\Y(\gl_{m|n})$, yielding the corresponding binary property in the super case, which was also pointed out in \cite{LM 2021}. As a consequence, we get  necessary and sufficient condition for the irreducibility of the $\Y(\gl_{m|n})$-module
\[
L_{h_1}(\la^{(1)})\ot L_{h_2}(\la^{(2)})\ot\cdots\ot L_{h_k}(\la^{(k)}),
\]
where each $\la^{(i)}$ is covariant for $1\le i\le k$. \ere

We have

\bco[Theorem B]
Let $\la_i$ be a covariant weight and  $h_i\in \BC$ for $1\leq i\leq k$.
Then the $\Y(\gl_{m|n})$-module
\[
L_{h_1}(\la^{(1)})\ot L_{h_2}(\la^{(2)})\ot\cdots\ot L_{h_k}(\la^{(k)}),
\]
is simple if and only if  $\la_i+h_i$ and $\la_j+h_j$ are super non-crossing for $1\leq i< j\leq k$.
\eco

\bigskip
\centerline{\bf Acknowledgments}
\medskip
V. Futorny is partially supported by  the National Natural Science Foundation of China (Grants 12350710787 and 12350710178). J.Zhang is partially supported by the National Natural Science Foundation of China (Grant 12571026).

\bibliographystyle{amsalpha}

\end{document}